\documentclass[11pt]{llncs}
\usepackage[a4paper,hmargin=2.5cm,vmargin=3cm]{geometry}

\usepackage{tikz}
\usetikzlibrary{arrows}

\usepackage[hyphens]{url}
\usepackage{amssymb,amsmath}
\usepackage{booktabs}
\usepackage{enumerate}

\newcommand{\F}{{\mathbb F}}
\newcommand{\K}{{\mathbb K}}
\newcommand{\N}{{\mathbb N}}

\DeclareMathOperator{\Gal}{Gal}
\DeclareMathOperator{\Gl}{Gl}

\begin{document}

\title{Breaking `128-bit Secure' Supersingular Binary Curves\thanks{The second author acknowledges the support 
of the Swiss National Science Foundation, via grant numbers 206021-128727 and 200020-132160,
while the third author acknowledges the support of the Irish Research Council, grant number ELEVATEPD/2013/82.}}
\subtitle{(or how to solve discrete logarithms in ${\mathbb F}_{2^{4 \cdot 1223}}$ and ${\mathbb F}_{2^{12 \cdot 367}}$)}

\author{Robert Granger\inst{1} \and Thorsten Kleinjung\inst{1} \and Jens Zumbr\"agel\inst{2}}

\institute{Laboratory for Cryptologic Algorithms, EPFL, Switzerland\\
\and
Institute of Algebra, TU Dresden, Germany\\
\email{robbiegranger@gmail.com, thorsten.kleinjung@epfl.ch, jens.zumbragel@ucd.ie}}

\maketitle

\begin{abstract}
In late 2012 and early 2013 the discrete logarithm problem (DLP) in finite fields of small characteristic underwent a dramatic series of breakthroughs,
culminating in a heuristic quasi-polynomial time algorithm, due to Barbulescu, Gaudry, Joux and Thom\'e.
Using these developments, Adj, Menezes, Oliveira and Rodr\'iguez-Henr\'iquez analysed the concrete security of the DLP, as it arises 
from pairings on (the Jacobians of) various genus one and two supersingular curves in the literature, which were originally thought to be $128$-bit 
secure. In particular, they suggested that the new algorithms have no impact on the security of a genus one curve over ${\mathbb F}_{2^{1223}}$, and 
reduce the security of a genus two curve over ${\mathbb F}_{2^{367}}$ to $94.6$ bits.
In this paper we propose a new field representation and efficient general descent principles which together make the new 
techniques far more practical. Indeed, at the `128-bit security level' our analysis shows that the aforementioned genus one curve
has approximately $59$ bits of security, and we report a total break of the genus two curve.


\end{abstract}

\keywords{Discrete logarithm problem, finite fields, supersingular binary curves, pairings}


\section{Introduction}\label{sec:intro}

The role of small characteristic supersingular curves in cryptography has been a varied and an interesting one. 
Having been eschewed by the cryptographic community for succumbing spectacularly to the subexponential 
MOV attack in 1993~\cite{MOV}, which maps the DLP from an elliptic curve (or more generally, the Jacobian of a higher genus curve) to 
the DLP in a small degree extension of the base field of the curve, they made a remarkable comeback with the advent of pairing-based 
cryptography in 2001~\cite{sakai,tripartite,bonehfranklin}. In particular, for the latter it was reasoned that the existence of a subexponential 
attack on the DLP does not {\em ipso facto} warrant their complete exclusion; rather, provided that the 
finite field DLP into which the elliptic curve DLP embeds is sufficiently hard, this state of affairs would be acceptable.

Neglecting the possible existence of native attacks arising from the supersingularity of these curves, much research effort  
has been expended in making instantiations of the required cryptographic operations on such curves as efficient as 
possible~\cite{barreto,stevenTate,duursmalee,granger04,granger05,ssab,1223,multicore,Chatter,ghosh,aranha,hasan}, to name but a few,
with the associated security levels having been estimated using Coppersmith's algorithm from 1984~\cite{coppersmith,Lenstra}.
Alas, a series of dramatic breakthrough results for the DLP in finite fields of small characteristic have potentially rendered all of these 
efforts in vain.

The first of these results was due to Joux, in December 2012, and consisted of a more efficient method --- dubbed `pinpointing' --- to obtain relations 
between factor base elements~\cite{Joux13a}.
For medium-sized base fields, this technique has heuristic complexity as low as 
$L(1/3,2^{1/3}) \approx L(1/3,1.260)$\footnote{The original paper states a complexity of 
$L(1/3,(8/9)^{1/3}) \approx L(1/3,0.961)$; however, on foot of recent communications the constant should be as stated.},
where as usual $L(\alpha,c) = L_Q(\alpha,c) = \text{exp}((c + o(1)) (\log{Q})^{\alpha} (\log{\log{Q} })^{1- \alpha})$, with $Q$ the 
cardinality of the field. This improved upon the previous best complexity of $L(1/3,3^{1/3}) \approx L(1/3,1.442)$ due to Joux and Lercier~\cite{JL06}.
Using this technique Joux solved example DLPs in fields of bitlength $1175$ and $1425$, both with prime base fields.

Then in February 2013, G\"olo\u{g}lu, Granger, McGuire and Zumbr\"agel used a specialisation of the Joux-Lercier doubly-rational 
function field sieve (FFS) variant~\cite{JL06}, in order to exploit a well-known family of `splitting polynomials', i.e., polynomials which
split completely over the base field~\cite{GGMZ13a}. For fields of the form $\F_{q^{kn}}$ with $k \ge 3$ fixed ($k =2$ is even simpler)
and $n \approx dq$ for a fixed integer $ d \ge 1$, 
they showed that for binary (and more generally small characteristic) fields,
relation generation for degree one elements runs in heuristic {\em polynomial time},
as does finding the logarithms of degree two elements (if $q^k$ can be written as $q'^{k'}$ for $k' \ge 4$), once degree one 
logarithms are known. For medium-sized base fields 
of small characteristic a heuristic complexity as low as $L(1/3,(4/9)^{1/3}) \approx L(1/3,0.763)$ was attained; this approach was 
demonstrated via the solution of example DLPs in the fields $\F_{2^{1971}}$~\cite{1971Ann} and $\F_{2^{3164}}$.

After the initial publication of~\cite{GGMZ13a}, Joux released a preprint~\cite{Joux13b}
detailing an algorithm for solving the discrete logarithm problem for fields of the form 
$\F_{q^{2n}}$, with $n \le q + d$ for some very small $d$, which was used to solve a DLP in $\F_{2^{1778}}$~\cite{1778Ann} 
and later in $\F_{2^{4080}}$~\cite{4080Ann}.
For $n \approx q$ this algorithm has heuristic complexity $L(1/4 + o(1),c)$ for some undetermined $c$, and also has a heuristic polynomial 
time relation generation method, similar in principle to that in~\cite{GGMZ13a}.
While the degree two element elimination method in~\cite{GGMZ13a} is arguably superior -- since elements can be eliminated on the fly --
for other small degrees Joux's elimination method is faster, resulting in the stated complexity.

In April 2013 G\"olo\u{g}lu {\em et al.} combined their approach with Joux's to solve an example DLP in the field $\F_{2^{6120}}$~\cite{6120Ann} 
and later demonstrated that Joux's algorithm can be tweaked to have heuristic complexity $L(1/4,c)$~\cite{GGMZ13b}, where 
$c$ can be as low as $(\omega/8)^{1/4}$~\cite{SACtalk}, with $\omega$ the linear algebra constant, i.e., the exponent of matrix multiplication.
Then in May 2013, Joux announced the solution of a DLP in the field $\F_{2^{6168}}$~\cite{6168Ann}.

Most recently, in June 2013, Barbulescu, Gaudry, Joux and Thom\'e announced a {\em quasi-polynomial time} for solving the DLP~\cite{barb}, 
for fields $\F_{q^{kn}}$ with $k \ge 2$ fixed and $n \le q + d$ with $d$ very small, which for $n \approx q$ has heuristic complexity
\begin{equation}\label{complexity}
(\log q^{kn})^{O(\log{\log{q^{kn}}})}.
\end{equation}
Since~(\ref{complexity}) is smaller than $L(\alpha,c)$ for any $\alpha > 0$, it is asymptotically the most efficient algorithm known for 
solving the DLP in finite fields of small characteristic, which can always be embedded into a field of the required form.
Interestingly, the algorithmic ingredients and analysis of this algorithm are much simpler than for Joux's $L(1/4 + o(1),c)$ algorithm.

Taken all together, one would expect the above developments to have a substantial impact on the security of small characteristic parameters appearing
in the pairing-based cryptography literature. However, all of the record DLP computations mentioned above used Kummer or twisted Kummer 
extensions (those with $n$ dividing $q^k \mp 1$), which allow for a reduction in the 
size of the factor base by a factor of $kn$ and make the descent phase for individual logarithms relatively easy. While such parameters are 
preferable for setting records (most recently in $\F_{2^{9234}}$~\cite{9234Ann}), none of the parameters featured in the 
literature are of this form, and so it is not {\em a priori} clear whether the new techniques weaken existing pairing-based protocol parameters.

A recent paper by Adj, Menezes, Oliveira and Rodr\'iguez-Henr\'iquez has begun to address this very issue~\cite{AMOR}. Using the time required to
compute a single multiplication modulo the cardinality of the relevant prime order subgroup as their basic unit of time, which we denote by $M_r$,
they showed that the DLP in the field $\F_{3^{6 \cdot 509}}$ costs at most $2^{73.7}$ $M_r$. One can arguably interpret this result to mean that 
this field has $73.7$ bits of security\footnote{The notion of bit security is quite fuzzy; for the elliptic curve DLP it is usually intended to mean the 
logarithm to the base $2$ of the expected number of group operations, however for the finite field DLP different authors have used different units, 
perhaps because the cost of various constituent algorithms must be amortised into a single cost measure. In this work we time everything in seconds, 
while to achieve a comparison with~\cite{AMOR} we convert to $M_r$.}.
This significantly reduces the intended security level of $128$ bits (or $111$ bits as 
estimated by Shinohara {{\em et al.}}~\cite{japan}, or $102.7$ bits for the Joux-Lercier FFS variant with pinpointing, as estimated in~\cite{AMOR}). 
An interesting feature of their 
analysis is that during the descent phase, some elimination steps are performed using the method
from the quasi-polynomial time algorithm of Barbulescu {\em et al.}, when one might have expected these steps to only come into play at much higher 
bitlengths, due to the high arity of the arising descent nodes.
 
In the context of binary fields, Adj {\em et al.} considered in detail the DLP in the field $\F_{2^{12 \cdot 367}}$,
which arises via a pairing from the DLP on the Jacobian of a supersingular genus two curve over $\F_{2^{367}}$, first proposed in~\cite{aranha},  
with embedding degree $12$. 
Using all of the available techniques they provided an upper bound of $2^{94.6}$ $M_r$ for the cost of breaking the DLP in the embedding field,
which is some way below the intended $128$-bit security level. In their conclusion Adj {\em et al.} also suggest that a commonly implemented 
genus one supersingular curve over $\F_{2^{1223}}$ with embedding degree~$4$~\cite{1223,multicore,Chatter,ghosh,hasan}, is not weakened by the 
new algorithmic advances, i.e., its security remains very close to $128$ bits.

In this work we show that the above security estimates were incredibly optimistic. Our techniques and results are summarised as follows.

\begin{itemize}
\item \textbf{Field representation:} We introduce a new field representation that can have a profound effect on the 
resulting complexity of the new algorithms. In particular it permits the use of a smaller $q$ than before, which not only speeds up the computation 
of factor base logarithms, but also the descent (both classical and new).
\vspace{2mm}
\item \textbf{Exploit subfield membership:} During the descent phase we apply a {\em principle of parsimony}, by which one 
should always try to eliminate an element in the target field, and only when this is not possible should one embed it into an extension field. 
So although the very small degree logarithms may be computed over a larger field, the descent cost is {\em greatly reduced} relative to 
solving a DLP in the larger field.
 \vspace{2mm}
\item \textbf{Further descent tricks:} The above principle also means that elements can automatically be rewritten in terms of elements of smaller degree, 
via factorisation over a larger field, and that elements can be eliminated via Joux's Gr\"obner basis computation method~\cite{Joux13b} 
with $k=1$, rather than $k > 1$, which increases its degree of applicability. 
\vspace{2mm}
\item \textbf{`128-bit secure' genus one DLP:} We show that the DLP in $\F_{2^{4 \cdot 1223}}$ can be solved in approximately $2^{40}\ \text{s}$,
or $2^{59}$ $M_r$, with $r$ a $1221$-bit prime.
\vspace{2mm}
\item \textbf{`128-bit secure' genus two DLP:} We report a total break of the DLP in $\F_{2^{12 \cdot 367}}$ (announced in~\cite{4404Ann}), 
which took about $52240$ core-hours.
\vspace{2mm}
\item $\mathbf{L(1/4,c)}$ \textbf{technique only:} Interestingly, using our approach the elimination steps \`a la  
Barbulesu {\em et al.}~\cite{barb} were not necessary for the above estimate and break.
\end{itemize}

The rest of the paper is organised as follows. In \S\ref{sec:setup} we describe our field representation and our target fields.
In \S\ref{sec:relgen} we present the corresponding polynomial time relation generation method for degree one elements and 
degree two elements (although we do not need the latter for the fields targeted in the present paper), 
as well as how to apply Joux's small degree elimination method~\cite{Joux13b} with the new representation.
We then apply these and other techniques to $\F_{2^{4 \cdot 1223}}$ in \S\ref{sec:genus1}  and to $\F_{2^{12 \cdot 367}}$ in \S\ref{sec:genus2} .
Finally, we conclude in \S\ref{sec:conclude}.






\section{Field Representation and Target Fields}\label{sec:setup}

In this section we introduce our new field representation and the fields whose DLP security we will address. This representation, as well as some 
preliminary security estimates, were initially presented in~\cite{ECCtalk}.

\subsection{Field Representation}\label{sec:fieldrep}

Although we focus on binary fields in this paper, for the purposes of generality, in this section we allow for extension fields of 
arbitrary characteristic. Hence let $q = p^l$ for some prime $p$, and let $\K = \F_{q^{kn}}$ be the field under consideration, with $k \ge 1$. 

We choose a positive integer $d_h$ such that $n \le q d_h  + 1$, and then choose $h_0,h_1 \in \F_{q^k}[X]$ with 
$\max \{\text{deg}(h_0),\text{deg}(h_1)\} = d_h$ such that 
\begin{equation}\label{fieldrep}
h_1(X^q) X - h_0(X^q) \equiv 0 \pmod{I(X)},
\end{equation}
where $I(X)$ is an irreducible degree $n$ polynomial in $\F_{q^k}[X]$. Then $\K = \F_{q^k}[X]/(I(X))$. 
Denoting by $x$ a root of $I(X)$, we introduce the auxiliary variable $y = x^q$, so that one has two isomorphic representations
of $\K$, namely $\F_{q^k}(x)$ and $\F_{q^k}(y)$, with $\sigma: \F_{q^k}(y) \rightarrow \F_{q^k}(x): y \mapsto x^q$.
To establish the inverse isomorphism, note that by~(\ref{fieldrep}) in $\K$ we have $h_1(y)x - h_0(y) = 0$, and hence 
$\sigma^{-1}: \F_{q^k}(x) \rightarrow \F_{q^k}(y): x \mapsto h_0(y)/h_1(y)$.

The knowledgeable reader will have observed that our representation is a synthesis of two other useful representations: the one used by 
Joux~\cite{Joux13b}, in which one searches for a degree $n$ factor $I(X)$ of $h_1(X)X^q - h_0(X)$;
and the one used by G\"olo\u{g}lu {\em et al.}~\cite{GGMZ13a,GGMZ13b}, in which one searches for a degree $n$ factor $I(X)$ of $X - h_0(X^q)$.
The problem with the former is that it constrains $n$ to 
be approximately $q$. The problem with the latter is that the polynomial $X - h_0(X^q)$ is insufficiently general to represent all degrees $n$ 
up to $q d_h$. By changing the coefficient of $X$ in the latter from $1$ to $h_1(X^q)$, we greatly increase the probability of overcoming the second
problem, thus combining the higher degree coverage of Joux's representation with the higher degree possibilities 
of~\cite{GGMZ13a,GGMZ13b}.

The {\em raison d'\^etre} of using this representation rather than Joux's representation is that for a given $n$, by choosing $d_h > 1$, one
may use a smaller $q$. So why is this useful? Well, since the complexity of the new descent methods is typically a function of $q$, then 
subject to the satisfaction of certain constraints, one may use a smaller $q$, thus reducing the complexity of solving the DLP.
This observation was our motivation for choosing field representations of the above form.

Another advantage of having an $h_1$ coefficient (which also applies to Joux's representation) is that it increases the chance of there being 
a suitable $(h_1,h_0)$ pair with coefficients defined over a proper subfield of $\F_{q^k}$, which then permits one to apply the factor base reduction technique 
of~\cite{JL06}, see~\S\ref{sec:genus1} and~\S\ref{sec:genus2}.


\subsection{Target Fields}\label{sec:target}

For $i \in \{0,1\}$ let $E_i/\F_{2^p}: Y^2 + Y = X^3 + X + i$. These elliptic curves are supersingular and can have prime or nearly prime order only 
for $p$ prime, and have embedding degree $4$~\cite{stevenSS,barreto,stevenTate}. 
We focus on the curve 
\begin{equation}\label{curve1}
E_{0} / \F_{2^{1223}}: Y^2 + Y = X^3 + X,
\end{equation}
which has a prime order subgroup of cardinality $r_{1} = (2^{1223} + 2^{612} +1)/5$, of bitlength $1221$.
This curve was initially proposed for $128$-bit secure protocols~\cite{1223} and has enjoyed several optimised 
implementations~\cite{multicore,Chatter,hasan,ghosh}. 
Many smaller $p$ have also been proposed in the literature (see~\cite{ssab,stevenSS}, for instance), and are clearly weaker.


For $i \in \{0,1\}$ let $H_i/\F_{2^p}: Y^2 + Y = X^5 + X^3 + i$. These genus two hyperelliptic curves are supersingular and can have a nearly 
prime order Jacobian only for $p$ prime (note that $13$ is always a factor of $\#\text{Jac}_{H_0}(\F_{2^p})$, since $\#\text{Jac}_{H_0}(\F_2) = 13$), 
and have embedding degree $12$~\cite{ssab,stevenSS}. 
We focus on the curve
\begin{equation}\label{curve2}
H_0 / \F_{2^{367}}: Y^2 + Y = X^5 + X^3,
\end{equation}
with $\#\text{Jac}_H(\F_{2^{367}}) = 13\cdot 7170258097 \cdot r_{2}$, and $r_2 = (2^{734} + 2^{551} + 2^{367} + 2^{184} + 1)/(13 \cdot 7170258097)$ is 
a $698$-bit prime,
since this was proposed for $128$-bit secure protocols~\cite{aranha}, and whose security was analysed in depth by Adj {\em et al.} in~\cite{AMOR}.



\section{Computing the Logarithms of Small Degree Elements}\label{sec:relgen}

In this section we adapt the polynomial time relation generation method from~\cite{GGMZ13a} and 
Joux's small degree elimination method~\cite{Joux13b} to the new field representation as detailed in~\S\ref{sec:fieldrep}.
Note that henceforth, we shall refer to elements of $\F_{q^{kn}} = \F_{q^k}[X]/(I(X))$ as field elements or as polynomials, as appropriate,
and thus use $x$ and $X$ (and $y$ and $Y$) interchangeably. We therefore freely apply polynomial ring concepts,
such as degree, factorisation and smoothness, to field elements.

In order to compute discrete logarithms in our target fields we apply
the usual index calculus method.  It consists of a precomputation
phase in which by means of (sparse) linear algebra techniques one
obtains the logarithms of the factor base elements, which will consist
of the low degree irreducible polynomials.  Afterwards, in the
individual logarithm phase, one applies procedures to recursively rewrite
each element as a product of elements of smaller degree, in this
way building up a {\em descent} tree, which has the target element as its root and factor base elements
as its leaves.  This proceeds in several stages,
starting with a continued fraction descent of the target element,
followed by a special-$Q$ lattice descent (referred to as degree-balanced classical
descent, see~\cite{GGMZ13a}), and finally using Joux's Gr\"obner basis descent~\cite{Joux13b} for the
lower degree elements. Details of the continued fraction and classical descent steps are 
given in \S\ref{sec:genus1}, while in this section we provide details of how to find the logarithms 
of elements of small degree.

We now describe how the logarithms of degree one and two elements (when needed) are to be computed. 
We use the relation generation method from~\cite{GGMZ13a}, rather than Joux's method~\cite{Joux13b}, since it automatically 
avoids duplicate relations. 
For $k \ge 2$ we first precompute the set $\mathcal{S}_k$, where
\[
\mathcal{S}_k = \{ (a,b,c) \in (\F_{q^k})^3 \mid X^{q+1} + a X^q + bX + c \ \ \text{splits completely over} \  \F_{q^k} \}.
\]
For $k=2$, this set of triples is parameterised by $(a,a^q, \F_{q} \ni c \ne a^{q+1})$, of which there are precisely $q^3 - q^2$ elements.
For $k \ge 3$, $\mathcal{S}_k$ can also be computed very efficiently, as follows. Assuming $ c \ne ab$ and 
$b \ne a^q$, the polynomial $X^{q+1} + a X^q + bX + c$ may be transformed (up to a scalar factor) into the polynomial
$f_B(\overline{X}) = \overline{X}^{q+1} + B\overline{X} + B$, where $B = \frac{(b - a^{q})^{q+1}}{(c - ab)^{q}}$, 
and $X = \frac{c - ab}{b - a^{q}} \overline{X} - a$. The set $\mathcal{L}$ of $B \in \F_{q^k}$ for which $f_B$ splits completely over $\F_{q^k}$ 
can be computed by simply testing for each such $B$ whether this occurs, and there are precisely $(q^{k-1}-1)/(q^2 - 1)$ such $B$ if $k$ is odd, and 
$(q^{k-1}-q)/(q^2 - 1)$ such $B$ if $k$ is even~\cite{Bluher}. Then for any $(a,b)$ such that $b \ne a^q$ and for each $B \in \mathcal{L}$, 
we compute via $B = \frac{(b - a^{q})^{q+1}}{(c - ab)^{q}}$ the corresponding (unique) $c \in \F_{q^{k}}$, which thus ensures that 
$(a,b,c) \in \mathcal{S}_k$. Note that in all cases we have $|\mathcal{S}_k| \approx q^{3k-3}$.

\subsection{Degree $1$ Logarithms}\label{sec:degree1}

We define the factor base $\mathcal{B}_1$ to be the set of linear elements in $x$, i.e., $\mathcal{B}_1 = \{ x - a \mid a \in \F_{q^k}\}$.
Observe that the elements linear in $y$ are each expressible in $\mathcal{B}_1$, since $(y - a) = (x - a^{1/q})^q$.

As in~\cite{JL06,GGMZ13a,GGMZ13b}, the basic idea is to consider elements of the form $xy + ay + bx + c$ with $(a,b,c) \in \mathcal{S}_k$. The above
two field isomorphisms induce the following equality in $\K$:
\begin{equation}\label{relation1}
x^{q+1} + a x^q + bx + c = \frac{1}{h_1(y)} \big(yh_0(y) + ay h_1(y) + bh_0(y) + ch_1(y) \big).
\end{equation}
When the r.h.s. of~(\ref{relation1}) also splits completely over $\F_{q^k}$, one obtains a relation between elements of $\mathcal{B}_1$ 
and the logarithm of $h_1(y)$. One can either adjoin $h_1(y)$ to the factor base, or simply use an $h_1(y)$ which splits completely over $\F_{q^k}$.

We assume that for each $(a,b,c) \in \mathcal{S}_k$ that the r.h.s. of~(\ref{relation1}) -- which has degree $d_h+1$ -- splits completely over $\F_{q^k}$ with probability $1/(d_h+1)!$. Hence in order for there to be sufficiently many relations we require that 
\begin{equation}\label{cond1}
\frac{q^{3k - 3}}{(d_h+1)!} > q^k, \ \ \text{or equivalently} \ \ q^{2k-3} > (d_h + 1)!.
\end{equation}
When this holds, the expected cost of relation generation is $(d_h+1)! \cdot q^{k} \cdot S_{q^k}(1,d_h+1)$, where $S_{q^k}(m,n)$ denotes
the cost of testing whether a degree $n$ polynomial is $m$-smooth, i.e., has all of its irreducible factors of degree $\le m$, 
see Appendix~B. The cost of solving the resulting linear system using sparse linear algebra techniques is 
$O(q^{2k+1})$ arithmetic operations modulo the order $r$ subgroup in which one is working.

\subsection{Degree $2$ Logarithms}\label{sec:degree2}

For degree two logarithms, there are several options. The simplest is to apply the degree one method over a quadratic extension of $\F_{q^k}$,
but in general (without any factor base automorphisms) this will cost $O(q^{4k+1})$ modular arithmetic operations. 
If $k \ge 4$ then subject to a condition on $q$, $k$ and $d_h$,
it is possible to find the logarithms of irreducible degree two elements on the fly, using the techniques of~\cite{GGMZ13a,GGMZ13b}. 
In fact, for the DLP in $\F_{2^{12 \cdot 367}}$ we use both of these approaches, but for different base fields, see~\S\ref{sec:genus2}.

Although not used in the present paper, for completeness we include here the analogue in our field representation of Joux's 
approach~\cite{Joux13b}. Since this approach forms the basis of the higher degree elimination steps in the quasi-polynomial time 
algorithm of Barbulescu {\em et al.}, its analogue in our field representation should be clear.

We define $\mathcal{B}_{2,u}$ to be the set of irreducible elements of $\F_{q^k}[X]$ of the form $X^2 + uX + v$. For each $u \in \F_{q^k}$ 
one expects there to be about $q^k/2$ such elements\footnote{For binary fields there are precisely $q^k/2$ irreducibles, since
$X^2 + uX + v$ is irreducible if and only if $\text{Tr}_{\F_{q^k}/\F_2}(v/u^2) = 1$.}. 
As in~\cite{Joux13b}, for each $u \in \F_{q^k}$ we find the logarithms of all 
the elements of $\mathcal{B}_{2,u}$ simultaneously. 
To do so, consider~(\ref{relation1}) but with $x$ on the l.h.s. replaced with $Q = x^2 + ux$. Using the field isomorphisms 
we have that $Q^{q+1} + a Q^q + bQ + c$ is equal to
\begin{align*}\label{relation2}
 &(y^2 \!+\! u^qy)((h_0(y)/h_1(y))^2 \!\!+\! u(h_0(y)/h_1(y))) \!+\! a(y^2 \!\!+\! u^qy) \!+\! b((h_0(y)/h_1(y))^2 \!\!+\! u(h_0(y)/h_1(y))) \!+\! c\\
 &= \frac{1}{h_1(y)^2} \big( (y^2 \!\!+\! u^qy)(h_0(y)^2 \!\!+\! uh_0(y)h_1(y)) \!+\! a(y^2 \!\!+\! u^qy)h_1(y)^2 \!+\! b(h_0(y)^2 \!\!+\! uh_0(y)h_1(y)) \!+\! ch_1(y)^2\big).
\end{align*}
The degree of the r.h.s. is $2(d_h+1)$, and when it splits completely over $\F_{q^k}$ we have a relation between elements of $\mathcal{B}_{2,u}$ and degree one elements,
whose logarithms are presumed known, which we assume occurs with probability $1/(2(d_h+1))!$.
Hence in order for there to be sufficiently many relations we require that 
\begin{equation}\label{cond2}
\frac{q^{3k - 3}}{(2(d_h+1))!} > \frac{q^k}{2}, \ \ \text{or equivalently} \ \ q^{2k-3} > (2(d_h + 1))!/2.
\end{equation}
Observe that~(\ref{cond2}) implies~(\ref{cond1}). When this holds, the expected cost of relation generation is 
$(2(d_h+1))! \cdot q^{k} \cdot S_{q^k}(1,2(d_h+1))/2$.
The cost of solving the resulting linear system using sparse linear algebra techniques is again $O(q^{2k+1})$
modular arithmetic operations, where now both the number of variables and the average weight is halved relative to the degree one case. 
Since there are $q^k$ such $u$, the total expected cost of this stage is $O(q^{3k+1})$ modular arithmetic operations, which may of course 
be parallelised.

\subsection{Joux's Small Degree Elimination with the New Representation}\label{sec:smallelim}

As in~\cite{Joux13b}, let $Q$ be a degree $d_Q$ element to be eliminated, let 
$F(X) = \sum_{i = 0}^{d_F} f_i X^i, G(X) = \sum_{j = 0}^{d_G} g_j X^j \in \F_{q^k}[X]$
with $d_F+d_G+2 \ge d_Q$, and assume without loss of generality $d_F \ge d_G$.
Consider the following expression:
\begin{equation}\label{corerelation}
G(X) \prod_{\alpha \in \F_q} ( F(X) - \alpha \, G(X) ) = F(X)^q G(X) - F(X) G(X)^q
\end{equation}
The l.h.s. is $\max(d_F,d_G)$-smooth.
The r.h.s. can be expressed modulo $h_1(X^q)X-h_0(X^q)$
in terms of $Y=X^q$ as a quotient of polynomials of relatively low degree
by using
\[
F(X)^q = \sum_{i = 0}^{d_F} f_{i}^q Y^i, \ G(X)^q = \sum_{j = 0}^{d_G} g_{j}^q Y^j
\ \text{and} \ X \equiv \frac{h_0(Y)}{h_1(Y)}.
\]
Then the numerator of the r.h.s. becomes
\begin{eqnarray}\label{fieldrelation}
\bigg( \sum_{i = 0}^{d_F} f_{i}^q Y^i \bigg)
\bigg( \sum_{j = 0}^{d_G} g_{j}^q h_0(Y)^j h_1(Y)^{d_F-j} \bigg)
-
\bigg( \sum_{i = 0}^{d_F} f_{i}^q h_0(Y)^i h_1(Y)^{d_F-i} \bigg)
\bigg( \sum_{j = 0}^{d_G} g_{j}^q Y^j \bigg).
\end{eqnarray}

Setting~(\ref{fieldrelation}) to be $0$ modulo $Q(Y)$ gives
a system of $d_Q$ equations over $\F_{q^k}$ in the $d_F+d_G+2$ variables $f_0,\ldots,f_{d_F},g_0, \ldots,g_{d_G}$.
By choosing a basis for $\F_{q^{k}}$ over $\F_q$ and
expressing each of the $d_F+d_G+2$ variables $f_0,\ldots,f_{d_F},g_0, \ldots,g_{d_G}$ in this
basis, this system becomes a bilinear quadratic system\footnote{The bilinearity makes finding solutions to this system easier~\cite{span},
and is essential for the complexity analysis in~\cite{Joux13b} and its variant in~\cite{GGMZ13b}.} 
of $kd_Q$ equations in $(d_F+d_G+2)k$ variables.
To find solutions to this system, one can specialise $(d_F+d_G+2-d_Q)k$ of the variables in order to make the resulting system generically zero-dimensional
while keeping its bilinearity, and then compute the corresponding Gr\"obner basis, which may have no solution, or a small number of solutions.
For each solution, one checks whether~(\ref{fieldrelation}) divided by $Q(Y)$ is $(d_Q-1)$-smooth: if so then 
$Q$ has successfully been rewritten as a product of elements of smaller degree; if no solutions give a $(d_Q-1)$-smooth cofactor, 
then one begins again with another specialisation.

The degree of the cofactor of $Q(Y)$ is upper bounded by $d_F (1 + d_h) - d_Q$, so assuming that it behaves as a uniformly
chosen polynomial of such a degree one can calculate the probability $\rho$ that it is $(d_Q-1)$-smooth using standard combinatorial techniques.

Generally, in order for $Q$ to be eliminable by this method with good probability, the number of solutions to the initial bilinear system must 
be greater than $1/\rho$.
To estimate the number of solutions, consider the action of $\Gl_2(\F_{q^{k}})$ on the set of pairs $(F,G)$.
The subgroups $\Gl_2(\F_q)$ and $\F_{q^{k}}^{\times}$ (via diagonal embedding) both
act trivially on the set of relations, modulo multiplication by elements in $\F_{q^{k}}^{\times}$. 
Assuming that the set of $(F,G)$ quotiented out by the action of the compositum of these subgroups 
(which has cardinality $\approx q^{k+3}$), generates distinct relations,  
one must satisfy the condition 
\begin{equation}\label{GBelimprob}
q^{(d_F+d_G +1 - d_Q)k - 3} > 1/ \rho\ .
\end{equation}
Note that while~(\ref{GBelimprob}) is preferable for an easy descent, one may yet violate it and still 
successfully eliminate elements by using various tactics, as demonstrated in~\S\ref{sec:genus2}.


\section{Concrete Security Analysis of $\F_{2^{4\cdot 1223}}$}%
\label{sec:genus1}

In this section we focus on the DLP in the $1221$-bit prime order $r_1$ subgroup of $\F_{2^{4\cdot 1223}}^{\times}$, 
which arises from the MOV attack applied to the genus one supersingular curve~(\ref{curve1}).
By embedding $\F_{2^{4\cdot 1223}}$ into its degree two
extension $\F_{2^{8\cdot 1223}} = \F_{2^{9784}}$ we show that, after
a precomputation taking approximately $2^{40}\ \text{s}$,
individual discrete logarithms can be computed in less than $2^{34}\ \text{s}$.
 

\subsection{Setup}

We consider the field $\F_{2^{8\cdot 1223}} = \F_{q^n}$ with $q = 2^8$
and $n = 1223$ given by the irreducible factor of degree $n$ of
$h_1(X^q)X - h_0(X^q)$, with
\[ h_0 = X^5 + tX^4 + tX^3 + X^2 + tX + t \:, \quad
h_1 = X^5 + X^4 + X^3 + X^2 + X + t \:, \]
where $t$ is an element of $\F_{2^2} \setminus \F_2$.
Note that the field of definition of this representation is $\F_{2^2}$.

Since the target element is contained in the subfield 
$\F_{2^{4\cdot 1223}}$, we begin the classical descent over
$\F_{2^4}$, we switch to $\F_q = \F_{2^8}$, i.e., $k=1$, for the
Gr\"obner basis descent, and, as explained below, we work over
$\F_{q^k}$ with either $k=1$ or a few $k>1$ to obtain the logarithms of
all factor base elements.




\subsection{Linear Algebra Cost Estimate}

In this precomputation we obtain the logarithms of all elements of
degree at most four over $\F_q$.  Since the degree $1223$ extension is defined over 
$\F_{2^2}$ in our field representation, by the action of
the Galois group $\Gal(\F_q/\F_{2^2})$ on the factor base, the number of irreducible elements of degree $j$ 
whose logarithms are to be computed can be reduced to about $2^{8j}/(4j)$ for $j \in \{1,2,3,4\}$.

One way to obtain the logarithms of these elements is to carry out the degree 1 relation generation method from~\S\ref{sec:degree1},
together with the elementary observation that an irreducible polynomial of degree $k$ over $\F_q$ splits completely over $\F_{q^k}$.
First, computing degree one logarithms over $\F_{q^3}$ gives the logarithms of
irreducible elements of degrees one and three over $\F_q$. Similarly,
computing degree one logarithms over $\F_{q^4}$ gives the logarithms of
irreducible elements of degrees one, two, and four over $\F_q$.
The main computational cost consists in solving the latter system arising from $\F_{q^4}$, which has size $2^{28}$ and an
average row weight of~$256$.

However, we propose to reduce the cost of finding these logarithms by using $k = 1$ only, in the following easy way.
Consider~\S\ref{sec:smallelim}, and observe that for each polynomial pair $(F,G)$ of degree at most~$d$, one obtains 
a relation between elements of degree at most $d$ when the numerator of the r.h.s. is $d$-smooth (ignoring factors of $h_1$).
Note that we are not setting the r.h.s. numerator to be zero modulo $Q$ or computing any Gr\"obner bases.
Up to the action of $\Gl_2(\F_q)$ (which gives equivalent relations) there
are about $q^{2d-2}$ such polynomial pairs.  Hence, for $d\ge 3$ there
are more relations than elements if the smoothness probability of the
r.h.s. is sufficiently high.  Notice that $k=1$ implies that the r.h.s. is
divisible by $h_1(Y)Y-h_0(Y)$, thus increasing its smoothness
probability and resulting in enough relations for $d=3$ and for $d=4$.
After having solved the much smaller system for $d=3$ we know the
logarithms of all elements up to degree three, so that the average row
weight for the system for $d=4$ can be reduced to about $\frac14 \cdot 256 =
64$ (irreducible degree four polynomials on the l.h.s.).  As above the
size of this system is $2^{28}$.

The cost for generating the linear systems is negligible compared to
the linear algebra cost.
For estimating the latter cost we consider
Lanczos' algorithm to solve a sparse $N\times N$,
$N=2^{28}$, linear system with average row weight $W=64$.  As noted
in~\cite{Popovyan,GGMZ13b} this algorithm can be
implemented such that
\begin{equation}\label{eq:lanczos}
  N^2 \,(2\,W\,\text{ADD} + 2\,\text{SQR} + 3\,\text{MULMOD})
\end{equation} operations are used.
On our benchmark system, an AMD Opteron 6168 processor at $1.9\,$GHz,
using~\cite{gmp} our implementation of these operations took 62 ns, 467 ns and 1853 ns
for an ADD, a SQR and a MULMOD, respectively,
resulting in a linear algebra cost of $2^{40}\ \text{s}$.

As in~\cite{AMOR}, the above estimate ignores communication costs and other possible slowdowns which may arise in practice.
An alternative estimate can be obtained by considering a problem of a similar size over $\F_2$ and 
extrapolating from~\cite{RSA768}. This gives an estimated time of $2^{42}\ \text{s}$, or for newer hardware slightly less.
Note that this computation was carried out using the block Wiedemann algorithm~\cite{wiedemann}, which we recommend
in practice because it allows one to distribute the main part of the computation.
For the sake of a fair comparison with~\cite{AMOR} we use the former estimate of $2^{40}\ \text{s}$.


\subsection{Descent Cost Estimate}\label{1223descent}

We assume that the logarithms of elements up to degree four are known, and that computing these logarithms with a lookup table is free.

\subsubsection{Small Degree Descent.}

We have implemented the small degree descent of \S\ref{sec:smallelim} in Magma~\cite{magma}
V2.20-1, using Faugere's F4 algorithm~\cite{F4}. For each degree from $5$ to $15$, on the same AMD Opteron 
6168 processor we timed the Gr\"obner basis computation between $10^6$ and $100$ times, depending on the degree.
Then using a bottom-up recursive strategy we estimated the following average running times in seconds for a full logarithm computation, 
which we present to two significant figures:
\[ C[5,\ldots,15] = 
[ \: 0.038 \,,\, 2.	1 \,,\, 2.1 \,,\, 93 \,,\, 95 \,,\, 180 \,,\,
190 \,,\, 3200 \,,\, 3500 \,,\, 6300 \,,\, 11000 \: ] \:. \]

\subsubsection{Degree-Balanced Classical Descent.}

From now on, we make the conservative assumption that a degree $n$ polynomial which is 
$m$-smooth, is a product of $n/m$ degree $m$ polynomials. In practice the descent 
cost will be lower than this, however, the linear algebra cost is dominating, so this issue is 
inconsequential for our security estimate. The algorithms we used for smoothness testing are detailed in Appendix B.

For a classical descent step with degree balancing we consider
polynomials $P(X^{2^a},Y) \in \F_q[X,Y]$ for a suitably chosen integer
$0 \le a \le 8$.  It is advantageous to choose $P$ such that its
degree in one variable is one; let $d$ be the degree in the other
variable.  In the case $\deg_{X^{2^a}}(P)=1$, i.e.,
$P=v_1(Y)X^{2^a}+v_0(Y)$, $\deg v_i \le d$, this gives rise to the
relation
\[ 
L_v^{2^a} = \bigg( \frac{R_v}{h_1(X)^{2^a}} \bigg)^{2^8} 
\quad \text{where}  \quad
\begin{array}{l}
  L_v = \tilde{v}_1(X^{2^{8-a}})X + \tilde{v}_0(X^{2^{8-a}}) \:, \\
  R_v = v_1(X)h_0(X)^{2^a}+v_0(X)h_1(X)^{2^a}
\end{array} 
\]
in $\F_q[X]/(h_1(X^q)X-h_0(X^q))$
with $\deg L_v \le 2^{8-a}d+1$, $\deg R_v \le d+5 \cdot 2^a$, and
$\tilde{v}_i$ being $v_i$ with its coefficients powered by $2^{8-a}$, for $i=0,1$.
Similarly, in the case $\deg_{Y}(P)=1$, i.e., $P=w_1(X^{2^a})Y+w_0(X^{2^a})$,
$\deg w_i \le d$, we have the relation
\[
L_w^{2^a} = \bigg(\frac{R_w}{h_1(X)^{2^a d}}\bigg)^{2^{8}} 
\quad \text{where} \quad
\begin{array}{l}
L_w =  \tilde{w}_1(X)X^{2^{8-a}}+ \tilde{w}_0(X) \:, \\
R_w = h_1(X)^{2^a d} \big(w_1 \big( \big(\frac{h_0(X)}{h_1(X)} \big)^{2^a}\big)X
+ w_0 \big( \big( \frac{h_0(X)}{h_1(X)} \big)^{2^a} \big) \big)
\end{array}
\]
with $\deg L_w \le d+2^{8-a}$,  $\deg R_w \le 5 \cdot 2^a d+1$ and
again $\tilde{w}_i$ being $w_i$ with its coefficients powered by $2^{8-a}$, for $i=0,1$.

The polynomials $v_i$ (respectively $w_i$) are chosen in such
a way that either the l.h.s.~or the r.h.s.~is
divisible by a polynomial $Q(X)$ of degree $d_Q$.  Gaussian reduction
provides a lattice basis $(u_0,u_1),(u_{0}',u_{1}')$ such that the
polynomial pairs satisfying the divisibility condition above are given
by $ru_i+su_{i}'$ for $i=0,1$, where $r,s \in \F_q[X]$.  For nearly all
polynomials $Q$ it is possible to choose a lattice basis of
polynomials with degree $\approx d_Q/2$ which we will assume for
all $Q$ appearing in the analysis; extreme cases can be avoided by
look-ahead or backtracking techniques.  Notice that a polynomial $Q$
over $\F_{2^4} \subset \F_q$ can be rewritten as a product of
polynomials which are also over $\F_{2^4} $, by choosing the basis as well as
$r$ and $s$ to be over $\F_{2^4}$.  This will be done in all steps of the
classical descent. The polynomials $r$ and $s$
are chosen to be of degree four, resulting in $2^{36}$ possible pairs
(multiplying both by a common non-zero constant gives the same relation).

In the final step of the classical eliminations (from degree $26$ to $15$)
we relax the criterion that the l.h.s. and r.h.s. are $15$-smooth, 
allowing also irreducibles of even degree up to degree $30$, since these can each be split over
$\F_q$ into two polynomials of half the degree, thereby
increasing the smoothness probabilities.  Admittedly, if we follow our worst-case analysis stipulation that all
polynomials at this step have degree $26$, then one could immediately split each of them into two degree $13$ polynomials.
However, in practice one will encounter polynomials of all degrees $\le 26$ and we therefore carry out the analysis 
without using the splitting shortcut, which will still provide an overestimate of the cost of this step.

In the following we will state the logarithmic cost (in seconds) of a
classical descent step as $c_l + c_r + c_s$, where $2^{c_l}$ and $2^{c_r}$ denote
the number of trials to get the left hand side and the right hand side
$m$-smooth, and $2^{c_s}\ \text{s}$ is the time required for the corresponding smoothness
test.  See Table~\ref{smoothness} for the smoothness timings that we
benchmarked on the AMD Opteron 6168 processor.

\begin{itemize}

\item $\mathbf{d_Q = 26}$ \textbf{to} $\mathbf{m = 15}$\textbf{:}
We choose $\deg_{X^{2^a}}P=1$, $a=5$, $Q$ on the right, and we have $d=17$,
$(\deg(L_v), \deg(R_v)) = (137, 151)$,
and logarithmic cost $13.4 + 15.6 - 9.0$, hence
$2^{20.0}\ \text{s}$; the expected number of factors is $19.2$,
so the subsequent cost will be less than $2^{17.7}\ \text{s}$.
Note that, as explained above, we use the splitting shortcut for irreducibles of even degree up to $30$, 
resulting in the higher than expected smoothness probabilities.
\vspace{2mm}
\item $\mathbf{d_Q = 36}$ \textbf{to} $\mathbf{m = 26}$\textbf{:} 
We choose $\deg_{X^{2^a}}P=1$, $a=5$, $Q$ on the right, and we have $d=22$,
$(\deg(L_v), \deg(R_v)) = (177, 146)$,
and logarithmic cost $18.7 + 13.6 - 9.0$, hence
$2^{23.3}\ \text{s}$; the expected number of factors is $12.4$,
so the subsequent cost will be less than $2^{23.9}\ \text{s}$.
\vspace{2mm}
\item $\mathbf{d_Q = 94}$ \textbf{to} $\mathbf{m = 36}$\textbf{:}
We choose $\deg_{Y}P = 1$, $a = 0$, $Q$ on the left, and we
have $d=51$, $(\deg(L_w), \deg(R_w))  = (213, 256)$, and
logarithmic cost $15.0 + 20.3 - 7.5$, hence $2^{27.8}\ \text{s}$; the
expected number of factors is $13.0$, so the subsequent cost
will be less than $2^{28.4}\ \text{s}$.\medskip
\end{itemize}

\vspace{-3mm}

\begin{table}[h]
  \caption{Timings for testing a degree $n$ polynomial over $\F_{2^4}$
    for $m$-smoothness.}
  \begin{center}\label{smoothness}
    \begin{tabular}{c|c|c}
      $n$ & $m$ & time \\
      \midrule      
      $137$ & $30$ & 1.9\,ms \\
      $146$ & $26$ & 1.9\,ms \\
      $213$ & $36$ & 5.1\,ms \\
      ~$611$~ & ~$94$~ & ~94\,ms \\
      \bottomrule
    \end{tabular}
  \end{center}
\end{table}

\vspace{-5mm}

\subsubsection{Continued Fraction Descent.}

For the continued fraction descent we multiply the target element by random powers of the generator
and express the product as a ratio of two polynomials of degree at most $611$. For each such expression we test
if both the numerator and the denominator are $94$-smooth. 
The logarithmic cost here is $17.7 + 17.7 - 3.4$, hence the cost is $2^{32.0}\ \text{s}$.
The expected number of degree $94$ factors on both sides will be
$13$, so the subsequent cost will be less than $2^{32.8}\ \text{s}$.

\subsubsection{Total Descent Cost}
The cost for computing an individual logarithm is therefore upper-bounded by $2^{32.0}\ \text{s} + 2^{32.8}\ \text{s} < 2^{34}\ \text{s}$.

\subsection{Summary}
The main cost in our analysis is the linear algebra computation which takes about $2^{40}\ \text{s}$,
with the individual logarithm stage being considerably faster. 
In order to compare with the estimate in \cite{AMOR}, we write the main cost in terms of $M_r$ which gives
$2^{59}$ $M_r$, and thus an improvement by a factor of $2^{69}$.
Nevertheless, solving a system of cardinality $2^{28}$ is still a formidable 
challenge, but perhaps not so much for a well-funded adversary.
For completeness we note that if one wants to avoid a linear algebra step of this size, then 
one can work over different fields, e.g., with $q = 2^{10}$ and $k=2$, or $q = 2^{12}$ and $k=1$.
However, while this allows a partitioning of the linear algebra into smaller steps as described in~\S\ref{sec:degree2}
but at a slightly higher cost, the resulting descent cost is expected to be significantly higher.

\section{Solving the DLP in $\F_{2^{12\cdot 367}}$}\label{sec:genus2}

In this section we present the details of our solution of a DLP in the $698$-bit prime order $r_2$ subgroup of 
$\F_{2^{12\cdot 367}}^{\times} =\F_{2^{4404}}^{\times}$, which arises from the MOV attack applied to the Jacobian of the genus 
two supersingular curve~(\ref{curve2}). Magma verification code is provided in Appendix A.
Note that the prime order elliptic curve $E_1 / \F_{2^{367}} : Y^2 + Y = X^3 + X + 1$ with embedding degree 4
also embeds into $\F_{2^{4404}}$, so that logarithms on this curve could have easily been computed as well.

\subsection{Setup}

To compute the target logarithm, as stated in~\S\ref{sec:intro} we applied
a principle of parsimony, namely, we tried to solve all intermediate logarithms in $\F_{2^{12 \cdot 367}}$,
considered as a degree $367$ extension of $\F_{2^{12}}$, and only when this was not possible did we embed
elements into the extension field $\F_{2^{24 \cdot 367}}$ (by extending the base field to $\F_{2^{24}}$) and 
solve them there. 

All of the classical descent down to degree 8 was carried out over $\F_{2^{12 \cdot 367}}$, which we
formed as the compositum of the following two extension fields. 
We defined $\F_{2^{12}}$ using the irreducible polynomial $U^{12} + U^3 + 1$ over $\F_2$, and defined $\F_{2^{367}}$ over $\F_2$
using the degree $367$ irreducible factor of $h_1( X^{64} )X - h_0( X^{64} )$, where $h_1 = X^5 + X^3 + X + 1$, and $h_0 = X^6 + X^4 + X^2 + X + 1$. 
Let $u$ and $x$ be roots of the extension defining polynomials in $U$ and $X$ respectively, and let $c = (2^{4404}-1)/r_2$.
Then $g = x + u^7$ is a generator of $\F_{2^{4404}}^{\times}$ and $\bar{g} = g^c$  is a generator of the subgroup of order $r_2$.
As usual, our target element was chosen to be $\bar{x}_{\pi} = x_{\pi}^{c}$ where 
\[
x_{\pi} = \sum_{i = 0}^{4403} (\lfloor \pi \cdot 2^{i+1} \rfloor \bmod 2)  \cdot  u^{11-(i \bmod 12)} \cdot  x^{ \lfloor i / 12 \rfloor}.
\]

The remaining logarithms were computed using a combination of tactics, over $\F_{2^{12}}$ when possible, and over $\F_{2^{24}}$  when not. 
These fields were constructed as degree 2 and 4 extensions of $\F_{2^{6}}$, respectively.
To define $\F_{2^{6}}$ we used the irreducible polynomial $T^6 + T +1$.
We then defined $\F_{2^{12}}$ using the irreducible polynomial $V^2 + tV + 1$ over $\F_{2^{6}}$, and 
$\F_{2^{24}}$ using the irreducible polynomial $W^4 + W^3 + W^2 + t^3$ over $\F_{2^{6}}$.

\subsection{Degree 1 Logarithms}

It was not possible to find enough relations for degree 1 elements over $\F_{2^{12}}$,  so in accordance
with our stated principle, we extended the base field to $\F_{2^{24}}$ to compute the logarithms of all $2^{24}$ 
degree 1 elements. We used the polynomial time relation generation from~\S\ref{sec:degree1}, which took 47 hours. 
This relative sluggishness was due to the r.h.s. having degree $d_h + 1 = 7$, which must split over $\F_{2^{24}}$.
However, this was faster by a factor of $24$ than it would have been otherwise, thanks to $h_0$ and $h_1$ being
defined over $\F_{2}$. This allowed us to use the technique from~\cite{JL06} to reduce the size of the factor base
via the automorphism $(x + a) \mapsto  (x + a)^{2^{367}}$, which fixes $x$ but has order $24$ on all non-subfield elements 
of $\F_{2^{24}}$, since $367 \equiv 7 \bmod 24$ and $\text{gcd}(7,24) = 1$. This reduced the factor base size to $699252$ elements, which 
was solved in 4896 core hours on a 24 core cluster using Lanczos' algorithm, approximately $24^2$ times faster than if we had not used 
the automorphisms.

\subsection{Individual Logarithm}

We performed the standard continued fraction initial split followed by degree-balanced classical descent as in~\S\ref{1223descent}, 
using Magma~\cite{magma} and NTL~\cite{ntl}, to reduce the target element to an 8-smooth product in $641$ and 
$38224$ core hours respectively. The most interesting part of the descent was the elimination of the elements of degree up to $8$
over $\F_{2^{12}}$ into elements of degree one over $\F_{2^{24}}$, which we detail below.
This phase was completed using Magma and took a further $8432$ core hours. However, we think that the combined time of 
the classical and non-classical parts could be reduced significantly via a backwards-induction analysis of the elimination 
times of each degree.

\subsubsection{Small Degree Elimination}

As stated above we used several tactics to achieve these eliminations. The 
first was the splitting of an element of even degree over $\F_{2^{12}}$ into two elements of half the degree (which had the same logarithm 
modulo $r_2$) over the larger field.
This automatically provided the logarithms of all degree 2 elements over $\F_{2^{12}}$.
Similarly elements of degree 4 and 8 over $\F_{2^{12}}$ were rewritten as elements of degree 2 and 4 over $\F_{2^{24}}$, 
while we found that degree 6 elements were eliminable more efficiently by initially continuing the descent over $\F_{2^{12}}$, as with degree 
$5$ and $7$ elements.

The second tactic was the application of Joux's Gr\"obner basis elimination method from~\S\ref{sec:smallelim} to elements 
over $\F_{2^{12}}$, as well as elements over $\F_{2^{24}}$. However, in many cases condition~(\ref{GBelimprob}) was violated, in which
case we had to employ various recursive strategies in order to eliminate elements. In particular, elements of the same degree were allowed 
on the r.h.s. of relations, and we then attempted to eliminate these using the same (recursive) strategy. For degree $3$ elements over
$\F_{2^{12}}$, we even allowed degree $4$ elements to feature on the r.h.s. of relations, since these were eliminable via the factorisation into
degree $2$ elements over $\F_{2^{24}}$.

In Figure 1 we provide a flow chart for the elimination of elements of degree up to $8$ over $\F_{2^{12}}$, and 
for the supporting elimination of elements of degree up to $4$ over $\F_{2^{24}}$. Nearly all of the arrows in Figure 1 were necessary 
for these field parameters (the exceptions being that for degrees 4 and 8 over $\F_{2^{12}}$ we could have initially continued the 
descent along the bottom row, but this would have been slower). The reason this `non-linear' descent arises is due to $q$ being 
so small, and $d_H$ being relatively large, which increases the degree of the r.h.s. cofactors, thus decreasing the smoothness probability. 
Indeed these tactics were only borderline applicable for these parameters; if $h_0$ or $h_1$ had 
degree any larger than $6$ then not only would most of the descent have been much harder, but it seems that one would be forced to compute 
the logarithms of degree 2 elements over 
$\F_{2^{24}}$ using Joux's linear system method from~\S\ref{sec:degree2}, greatly increasing the required number of core hours.
As it was, we were able to eliminate degree 2 elements over $\F_{2^{24}}$ on the fly, as we describe explicitly below.

Finally, we note that our descent strategy is considerably faster than the alternative of embedding the DLP into $\F_{2^{24 \cdot 367}}$ 
and performing a full descent in this field, even with the elimination on the fly of degree 2 elements over $\F_{2^{24}}$, since much of the resulting 
computation would constitute superfluous effort for the task in hand.

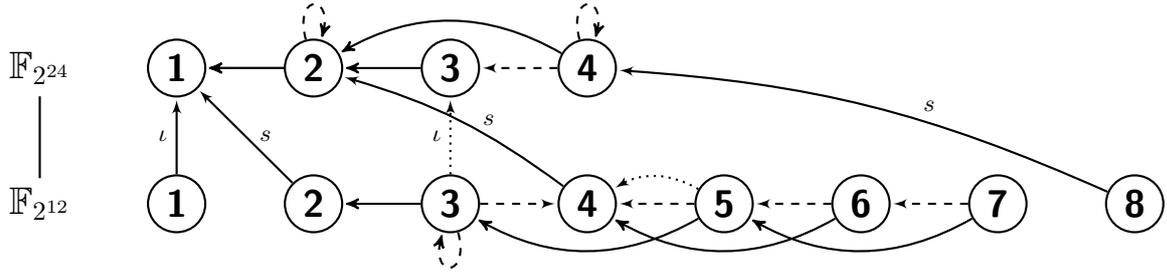
\begin{figure}[t]\label{GBpic}
\tikzstyle{line} = [draw, -latex']
\begin{tikzpicture}[->,>=stealth',shorten >=1pt,auto,node distance=1.8cm,
  thick,main node/.style={circle,draw,font=\sffamily\Large\bfseries}, scale=1, transform shape]

  \node[main node] (1) {1};
  \node[main node] (2) [right of=1] {2};
  \node[main node] (3) [right of=2] {3};
  \node[main node] (4) [right of=3] {4};
  \node[main node] (1a) [below of=1] {1};
  \node[main node] (2a) [below of=2] {2};
  \node[main node] (3a) [below of=3] {3};
  \node[main node] (4a) [below of=4] {4};
  \node[main node] (5a) [right of=4a] {5};
  \node[main node] (6a) [right of=5a] {6};
  \node[main node] (7a) [right of=6a] {7};
  \node[main node] (8a) [right of=7a] {8};
  \node (f24) [left of=1] {\Large{$\F_{2^{24}}$}};
  \node (f12) [left of=1a] {\Large{$\F_{2^{12}}$}};
   \path [line] (f12) edge [-] (f24);

   \path [line] (1a) edge node [left] {$\iota$} (1);
   \path [line,dotted] (3a) edge node [left] {$\iota$} (3);
   \path [line] (2a) edge node [right] {\hspace{0.7mm}$s$} (1);
  \path [line,dashed] (3a) edge [loop below] node {} (3a);
  \path [line,dashed] (3a) edge node {} (4a);
   \path [line] (4a) edge [bend right=10] node [right] {\hspace{2mm}$s$} (2);
  \path [line,dashed] (6a) edge node {} (5a);
  \path [line,dashed] (7a) edge node {} (6a);
   \path [line] (8a) edge [bend right=10] node [right] {\hspace{6mm}$s$} (4);
  \path [line,dotted] (5a) edge [bend right] node {} (4a);
  \path [line,dashed] (4) edge [loop above] node {} (4);
  \path [line,dashed] (4) edge node {} (3);

  \path [line,dashed] (5a) edge node {} (4a);
  \path [line,dashed] (2) edge [loop above] node {} (2);
 \path[every node/.style={font=\sffamily\small}]
       (2) edge node {} (1)
       (3) edge node {} (2)
       (2) edge node {} (1)
       (3a) edge node {} (2a)
       (6a) edge [bend left ] node {} (4a)
       (5a) edge [bend left ] node {} (3a)
       (4) edge [ bend right] node {} (2)
       (7a) edge [ bend left] node {} (5a);
\end{tikzpicture}
\caption{This diagram depicts the set of strategies employed to eliminate elements over $\F_{2^{12}}$ of degree up to $8$. 
The encircled numbers represent the degrees of elements over $\F_{2^{12}}$ on the bottom row, and over $\F_{2^{24}}$ on the top row.
The arrows indicate how an element of a given degree is rewritten as a product of elements of other degrees, possibly over the larger field.
Unadorned solid arrows indicate the maximum degree of elements obtained on the l.h.s. of the Gr\"obner basis elimination method; likewise dashed
arrows indicate the degrees of elements obtained on the r.h.s. of the Gr\"obner basis elimination method, when these are greater than those obtained 
on the l.h.s. Dotted arrows indicate a fall-back strategy when the initial strategy fails.
An $s$ indicates that the element is to be split over the larger field into two elements of half the degree.
An $\iota$ indicates that an element is promoted to the larger field. Finally, a loop indicates that one must use a recursive strategy in which further
instances of the elimination in question must be solved in order to eliminate the element in question.
}
\end{figure}

\subsubsection{Degree 2 Elimination over $\F_{2^{24}}$}

Let $Q(Y)$ be a degree two element which is to be eliminated, i.e., written as a product of degree one elements.
As in~\cite{GGMZ13a,GGMZ13b} we first precompute the set of $64$ elements $B \in \F_{2^{24}}$ such that the polynomial 
$f_B(X) = X^{65} + BX + B$ splits completely over $\F_{2^{24}}$
(in fact these $B$'s happen to be in $\F_{2^{12}}$, but this is not relevant to the method).
We then find a Gaussian-reduced basis of the lattice $L_{Q(Y)}$ defined by 
\[ 
L_{Q(Y)} = \{(w_0(Y),w_1(Y)) \in \F_{2^{24}}[Y]^2 : w_0(Y) \,h_0(Y) + w_1(Y)\, h_1(Y) \equiv 0 \pmod{Q(Y)}\} \:. 
\]
Such a basis has the form $(u_{0}, Y+u_{1}), (Y+v_{0}, v_{1})$, with $u_i,v_i \in \F_{2^{24}}$,
except in rare cases, see Remark 1.
For $s\in\F_{2^{24}}$ we obtain lattice elements $(w_0(Y), w_1(Y)) = (Y + v_{0} + su_{0}, sY + v_{1} + su_{1})$. 

Using the transformation detailed in \S\ref{sec:relgen}, for each $B \in \F_{2^{24}}$ such that $f_B$ splits completely over $\F_{2^{24}}$
we perform a Gr\"obner basis computation to find the set of $s \in \F_{2^{24}}$ that satisfy
\[
B = \frac{( s^{64} + u_{0}s + v_{0})^{65}}
{(u_{0}s^2 + (u_{1} + v_{0})s+ v_{1} )^{64}} \:, 
\]
by first expressing $s$ in a $\F_{2^{24}} / \F_{2^{6}}$ basis, which results in a quadratic system in $4$ variables.
This ensures that the l.h.s. splits completely over $\F_{2^{24}}$. 
For each such $s$ we check whether the r.h.s. cofactor of $Q(Y)$, which has degree $5$, is $1$-smooth.
If this occurs, we have successfully eliminated $Q(Y)$.

However, one expects on average just one $s$ per $B$, and so the probability of $Q(Y)$ being eliminated in this way is 
$1 - (1 - 1/5!)^{64} \approx 0.415$, which was borne out in practice to two decimal places.
Hence, we adopted a recursive strategy in which we stored all of the r.h.s. cofactors 
whose factorisation degrees had the form $(1,1,1,2)$ (denoted type 1), or $(1,2,2)$ (denoted type 2). Then for each type 1 
cofactor we checked to see if the degree 2 factor was eliminable by the above method. If none were eliminable
we stored every type 1 cofactor of each degree 2 irreducible occurring in the list of type 1 cofactors of $Q(Y)$. 
If none of these were eliminable (which occurred with probability just $0.003$), then we reverted to the type 2 
cofactors, and adopted the same strategy just specified for each of the degree 2 irreducible factors.
Overall, we expected our strategy to fail about once in every $6 \cdot 10^6$ such $Q(Y)$. This happened just once during our descent, 
and so we multiplied this $Q(Y)$ by a random linear polynomial over 
$\F_{2^{24}}$ and performed a degree 3 elimination, which 
necessitates an estimated 32 degree 2 polynomials being simultaneously eliminable by the above method, which 
thanks to the high probability of elimination, will very likely be successful for any linear multiplier.


\subsection{Summary}

Finally, after a total of approximately 52240 core hours (or $2^{48}$ $M_{r_2}$), 
we found that $\bar{x}_{\pi} = \bar{g}^{\text{log}}$, with $\text{log} =$
\begin{align*}
&40932089202142351640934477339007025637256140979451423541922853874473604
\\[-.5mm]
&39015351684721408233687689563902511062230980145272871017382542826764695
\\[-.5mm]
&59843114767895545475795766475848754227211594761182312814017076893242 \:.
\\[-.5mm]
\end{align*}

\begin{remark}
During the descent, we encountered several polynomials $Q(Y)$ that were apparently not eliminable via the Gr\"obner basis method.
We discovered that they were all factors of $h_1(Y) \cdot c + h_0(Y)$ for $c \in \F_{2^{12}}$ or $\F_{2^{24}}$, and hence 
$h_0(Y)/h_1(Y) \equiv c \pmod{Q(Y)}$.
This implies that~(\ref{fieldrelation}) is equal to $F(c)G^{(q)}(Y)+F^{(q)}(Y)G(c)$ modulo $Q(Y)$,
where $G^{(q)}$ denotes the Frobenius twisted $G$ and
similarly for $F^{(q)}$.
This cannot become $0$ modulo $Q(Y)$ if the degrees
of $F$ and $G$ are smaller than the degree of $Q$, unless $F$ and $G$ are both constants.
However, thanks to the field representation, finding the logarithm of these $Q(Y)$ turns out to be easy. In particular, if 
$h_1(Y) \cdot c + h_0(Y) = Q(Y) \cdot R(Y)$ then
$Q(Y) = h_1(Y) \cdot ((h_0/h_1)(Y) + c)/R(Y) = h_1(Y) \cdot (X + c) / R(Y)$, and thus modulo $r_2$ we have 
$\log(Q(y)) \equiv \log(x+c) - \log(R(y))$, since $\log(h_1(y)) \equiv 0$. Since $(x + c)$ is in the factor base, if we are able to compute 
the logarithm of $R(y)$, then we are done. In all the cases we encountered, the cofactor $R(y)$ was solvable by the above methods.
\end{remark}

\section{Conclusion}\label{sec:conclude}

We have introduced a new field representation and efficient descent principles which together make the recent DLP advances far more practical.
As example demonstrations, we have applied these techniques to two binary fields of central interest to pairing-based cryptography, namely 
$\F_{2^{4 \cdot 1223}}$ and $\F_{2^{12 \cdot 367}}$, which arise as the embedding fields of (the Jacobians of) a genus one and a genus two 
supersingular curve, respectively. When initially proposed, these fields were believed to be $128$-bit secure, and even in light of the recent 
DLP advances, were believed to be $128$-bit and $94.6$-bit secure. On the contrary, our analysis indicates that the former field has approximately 
$59$ bits of security and we have implemented a total break of the latter.


\bibliographystyle{plain}
\bibliography{crypto}


\section*{Appendix A}\label{sec:verify}

The following Magma script verifies the solution of the chosen DLP in the order $r_2$ subgroup of $\F_{2^{12 \cdot 367}}^{\times}$:
{
\small{
\begin{verbatim}
// Field setup
F2 := GF(2);
F2U<U> := PolynomialRing(F2);
F2_12<u> := ext< F2 | U^12 + U^3 + 1 >;
F2_12X<X> := PolynomialRing(F2_12);

modulus := (2^734 + 2^551 + 2^367 + 2^184 + 1) div (13 * 7170258097);
cofactor := (2^4404 - 1) div modulus;

h1 := X^5 + X^3 + X + 1;
h0 := X^6 + X^4 + X^2 + X + 1;
temp1 := Evaluate(h1, X^64) * X + Evaluate(h0, X^64);
temp2 := X^17 + X^15 + X^14 + X^13 + X^12 + X^11 + X^10 + X^6 + 1;
polyx := temp1 div temp2;
Fqx<x> := ext< F2_12 | polyx >;

// This is a generator for the entire multiplicative group of GF(2^4404).
g := x + u^7;

// Generate the target element.
pi := Pi(RealField(1500));
xpi := &+[ (Floor(pi * 2^(i+1)) mod 2) * u^(11-(i mod 12)) * x^(i div 12) : i in [0..4403]];

log := 4093208920214235164093447733900702563725614097945142354192285387447\
36043901535168472140823368768956390251106223098014527287101738254282676469\
559843114767895545475795766475848754227211594761182312814017076893242;

// If the following is true, then verification was successful
(g^cofactor)^log eq xpi^cofactor;
\end{verbatim}
}}

\section*{Appendix B}\label{sec:smoothness}

This section provides the algorithmic details of the smoothness testing function used in 
\S\ref{sec:genus1} and \S\ref{sec:genus2}.
Given a polynomial $f(X)$ of degree $n$ over $\F_q$, in order to test
its $m$-smoothness we compute
\[ t(X) \::=\: f'(X) \prod_{\lfloor m/2 \rfloor + 1}^m 
(X^{q^i} - X) \mod f \:. \]

Let $R$ be the quotient ring $\F_q[X] / \langle f \rangle$ (so that $R
\cong \F_q^n$ as vector spaces), and denote a residue class in $R$ by
$[a(X)]$.  A multiplication in $R$ can be computed using $2n^2$
$\F_q$-multiplications.
In order to obtain the above product our main task is to compute
$[X^{q^i}]$ for $i \in \{ \lfloor m/2 \rfloor + 1, \dots, m\}$, after
which we can compute $t(X)$ using $\lceil m/2 \rceil$
$R$-multiplications.


\subsection*{How to Compute a Power $[X^{p^{rs}}]$}

First we explain a method how to obtain a general power $[ X^{p^{rs}}
]$, where $p$ is the characteristic of $\F_q$.
We precompute $[ X^{p^r} ], [ X^{2p^r} ], \dots, [ X^{(n-1)p^r} ]$ by
consecutively multiplying by $[ X ]$ (i.e., shifting).  This requires
$(n-1)(p^r-1)$ shifts, each using $n$ $\F_q$-multiplications, so less
than $n^2 p^r$ $\F_q$-multiplications in total.

With this precomputation we then can compute $p^r$-powering in $R$,
i.e., one application of the map $\varphi:R\to R$,
$\alpha\to\alpha^{p^r}$, in the following way:
\[ \big[ \sum_{i=0}^{n-1} a_i X^i \big] ^ {p^r} =
\sum_{i=0}^{n-1} a_i^{p^r} \big[ X^{ip^r} \big] \]
This requires $n$ $p^r$-powering operations in $\F_q$ (which we ignore) and $n$
scalar multiplications in $R$, hence $n^2$ $\F_q$-multiplications.
Finally, we compute the powers $[ X^{p^{ri}} ]$ by repeatedly
applying the map~$\varphi$, i.e., $[ X^{p^{ri}} ] = \varphi^i([X]) =
\varphi^{i-1}([X^{p^r}])$, for $i\in\{2,\dots,s\}$, which requires
$(s-1)$ $p^r$-powerings in $R$.  Altogether we can compute $[
X^{p^{rs}} ]$ in less than $n^2(p^r + s)$ $\F_q$-operations.

\medskip

For an alternative method of computing $[ X^{p^r} ], [ X^{2p^r} ],
\dots, [ X^{(n-1)p^r} ]$ we assume that $[ X^{p^r} ]$ is already
known.  First, by multiplying by $[X]$ we obtain $[ X^{p^r+1} ], [
X^{p^r+2} ], \dots, [ X^{p^r+(n-1)} ]$ using $(n-1)$ shifts, hence
less than $n^2$ $\F_q$-multiplications.  With this we can compute a
multiplication by $X^{p^r}$, i.e.,
\[ \big[ \sum_{i=0}^{n-1} a_i X^i \big] \cdot [ X^{p^r} ] =
\sum_{i=0}^{n-1} a_i \big[ X^{p^r+i} \big] \:, \]
using $n^2$ $\F_q$-multiplications.  We apply this multiplication map
repeatedly in order to compute $[ X^{p^r} ], [ X^{2p^r} ], \dots, [
X^{(n-1)p^r} ]$; instead of $n^2 p^r$ $\F_q$-multiplications, 
this method requires $n^3$ $\F_q$-multiplications.


\subsection*{Computing the Powers $[X^{q^i}]$}

We outline two strategies to compute the powers $[X^{q^i}]$ for
$i\in\{1,\dots,m\}$.

\paragraph{Strategy 1}

Write $q = (p^r)^s = p^{rs}$.  As in the method outlined above
we do a precomputation in order to represent the $p^r$-powering
map in $R$.  We then apply this map repeatedly
in order to compute $[ X^{p^{rj}} ]$ for $j\in\{2,\dots,sm\}$,
and obtain this way the powers $[X^{q^i}] = [ X^{p^{rsi}} ]$.

This method requires about $n^2(p^r+sm)$ $\F_q$-multiplications.

\paragraph{Strategy 2}

First we compute $[X^q]$ by writing $q = p^{rs}$ and using the above
method, which requires $n^2(p^r+s)$ $\F_q$-multiplications.
With this we can use the alternative method outlined above
for precomputing the $q$-powering map in $R$; here
we let $s=1$, i.e., $q=p^r$.  We then apply this map
repeatedly to obtain the powers $[X^{q^i}]$.
This method requires about $n^2(p^r+s+n+m)$ $\F_q$-multiplications, and corresponds to the smoothness test in the Adj {\em et al.}
paper; but the version here has an improved running time
(the previous one was $n^2(2n+m+4\log q)$ $\F_q$-multiplications).


\subsection*{Examples}

In the case $q = 2^8$ the running time (in $\F_q$-multiplications) using Strategy~1 and $s = 2$ is $n^2(16+2m)$, while using Strategy~2 and $s = 4$ it is $n^2(8+n+m)$.  When $q = 2^4$ the running time using Strategy~1 and $s = 1$ is $n^2(16+m)$, and using Strategy~2 and $s = 2$ is $n^2(6+n+m)$.  Hence, for typical values of $n$ and $m$ we prefer and implement Strategy~1.  For example, if $q = 2^4$, $n = 611$, $m = 94$ (see \S4.3) we need $110 n^2$ $\F_q$-multiplications.

%
%
%


\begin{remark}
Recall that in either case, in order to obtain $t(X)$ and thus to
complete the smoothness test, we have to consider the final $\lceil
m/2 \rceil$ $R$-multiplications.  This requires an additional cost of
about $n^2m$ $\F_q$-multiplications.
\end{remark}

\end{document}